\documentclass[12pt,reqno,a4paper]{amsart}
\usepackage{amsmath,amssymb,amsfonts,amsthm}
\usepackage[mathscr]{eucal}

\author{H.~M.~Khudaverdian}
\author{Th.~Th. Voronov}

\address{G.~S.~Sahakian~Department~of~Theoretical ~Physics, Yerevan State
University, 1 A. Manoukian Street, 375049 Yerevan, Armenia}

\address{Department of Mathematics, University of Manchester Institute of Science and Technology
(UMIST), United Kingdom}

\email{theodore.voronov@umist.ac.uk, khudian@umist.ac.uk}

\title[Odd Laplace operators]{Geometry of differential operators, and odd Laplace operators}

\newtheorem{thm}{Theorem}

\newtheorem*{coro}{Corollary}

\theoremstyle{definition}

\renewcommand{\leq}{\leqslant}

 \DeclareMathOperator{\ord}{ord}
 \DeclareMathOperator{\Res}{Res}

\newcommand{\Act}{{\mathscr{A}}}

\DeclareMathOperator{\Vol}{Vol} 
 \DeclareMathOperator{\Vect}{Vect}
\DeclareMathOperator{\sub}{sub}

\newcommand{\lie}[1]{{\mathcal L}_{{#1}}}

\newcommand{\der}[2]{{\frac{\partial {#1}}{\partial {#2}}}}

\newcommand{\RR}{\mathbb R}

\newcommand{\p}{\partial}
\newcommand{\fun}{C^{\infty}}
\newcommand{\V}{{\mathfrak{V}}}

\def\e{\varepsilon}

\def\D{\Delta}

\newcommand{\g}{{\gamma}}

\renewcommand{\S}{{{S}}}
\newcommand{\lt}{\theta} 

\newcommand{\ft}{{\tilde f}}

\newcommand{\at}{{\tilde a}}
\newcommand{\bt}{{\tilde b}}

\newcommand{\ps}{{\boldsymbol{\psi}}}
\newcommand{\ch}{{\boldsymbol{\chi}}}

\begin{document}

\maketitle

Let  $\D$ be an arbitrary linear differential operator of the second
order acting on functions on a (super)manifold  $M$. In local
coordinates $\D=\frac{1}{2}\,S^{ab}\,\p_b\p_a+ T^a\,\p_a +R$. The
principal symbol of  $\D$ is the symmetric tensor field $S^{ab}$, or
the quadratic function $S=\frac{1}{2}\,S^{ab} p_bp_a$ on $T^*M$. The
principal symbol can be understood as a symmetric ``bracket'' on
functions: $\{f,g\}:=\D(fg)-(\D f)\,g-(-1)^{\e\ft}f\,(\D g)
+\D(1)\,fg$, where $\e=\tilde\D$ is the parity of the operator $\D$;
in coordinates $\{f,g\}=S^{ab}\p_bf\,\p_a g (-1)^{\at\ft}$. In the
following by a \textit{bracket} in a commutative algebra we mean an
arbitrary symmetric bi-derivation.  The problem is to describe   all
operators $\D$ with a given $S^{ab}$, or, which is the same,   all
operators generating a given bracket $\{f,g\}$. Without loss of
generality we set $R=\D(1):=1$ in the sequel. Initially we suppose
that the operators   act on scalar functions; operator pencils acting
on densities of arbitrary weights will naturally appear in the course
of study. Everything is applicable to supermanifolds as well as to
usual manifolds. For odd operators in the super case questions about
identities of the Jacobi type arise. The problem is closely related
with the geometry of the Batalin--Vilkovisky formalism in quantum
field theory (description of the ``generating operators'' for an odd
bracket).

The first non-trivial observation is that {H\"ormander}'s
subprincipal symbol $\sub\D=(\p_bS^{ba}(-1)^{\bt(\e+1)}-2T^a)p_a$ can
be interpreted as  an  ``upper connection'' in the bundle $\Vol M$.
Precisely, $\g^a=\p_bS^{ba}(-1)^{\bt(\e+1)}-2T^a$ has the
transformation law $\g^{a'}=\left(\g^a+S^{ab}\,\p_b\ln
J\right)\der{x^{a'}}{x^a}$, where $J=\frac{Dx'}{Dx}$ (the Jacobian),
and it specifies a ``contravariant derivative''
$\nabla^a\rho=(S^{ab}\p_b +\g^a)\rho$ on volume forms. The
coordinate-dependent Hamiltonian $\g=\sub\D=\g^a p_a$ plays the role
of a local connection form. If the matrix  $S^{ab}$ is invertible,
then we can lower the index $a$ and get a usual connection. (Let us
stress that  $\D$ acts on functions, and \textit{a priori} there is
no  extra structure on our manifold. The bundle $\Vol M$ and an upper
connection in it arise from the operator itself.) Thus, $\D$ is
defined by a set of data: a bracket on functions and an associated
upper connection in $\Vol M$.

Define the \textit{algebra of densities} $\V(M)$ as the algebra of
formal linear combinations of densities of arbitrary weights $w\in
\RR$. In  $\V(M)$ there is a unit $1$ and a natural invariant
scalar multiplication. The scalar product is given by the
formula: $\langle \ps,\ch\rangle=\int_M \Res
(t^{-2}\ps(x,t)\ch(x,t))\,Dx$.  We specify elements of  $\V(M)$
by generating functions $\ps(x,t)$ defined on a manifold $\hat
M$. One can classify the derivations of
$\V(M)$~\cite{tv:laplace2}. A {bracket}   of weight $0$ in
$\V(M)$ is specified by a tensor $(\hat S^{\hat a\hat b})=
    \bigl(\begin{smallmatrix}
       S^{ab} & t  \g^a   \\
      t  \g^a & t^2\lt   \\
    \end{smallmatrix}
    \bigr)$ on $\hat M$, where $S^{ab}$ gives a bracket on $M$,
$\g^a$ gives an upper connection in $\Vol M$ associated with
$S^{ab}$, and the term $\lt$ is a second order geometric object
(depending on  $S^{ab}$ and $\g^a$), similar to the Brans--Dicke
field in Kaluza--Klein type theories. A set of data $S^{ab}, \g^a,
\lt$ is equivalent (non-canonically) to a set consisting of $S^{ab}$,
a vector field and a scalar. Operators in the algebra $\V(M)$ are
written as operator pencils $\D_w$ acting on $w$-densities. A pencil
$\D_w$ is self-adjoint if $(\D_w)^*=\D_{1-w}$. A second order
differential operator in $\V(M)$ is represented by a quadratic pencil
$\D_w=\D_0+wA+w^2B$, where  $\D_0$ is a second order  operator on
functions,  $A$ and $B$ having orders $\leq 1$ and $0$.

\begin{thm}
For the algebra $\V(M)$ there is a one-to-one correspondence
between brackets and second order operators with the
self-adjoint\-ness condition. To a bracket with the matrix  $\hat
S^{\hat
a\hat b}$ corresponds a  ``canonical pencil''\\
 $\D_w=\frac{1}{2}\left(\S^{ab}\p_b\p_a+
    \left(\p_b\S^{ba}(-1)^{\bt(\e+1)}+(2 w -1)\g^a\right)\p_a
    +\right. \\
    \hfill \left.
     w\,\p_a\g^a(-1)^{\at(\e+1)}  +
         w(w -1)\,\lt \right).
         $
\end{thm}
The corresponding operator in the algebra $\V(M)$ is the Laplacian
constructed from   $\hat S^{\hat a\hat b}$ and the canonical
divergence on $\hat M$ ~\cite{tv:laplace2}. There is a unique
canonical pencil passing through an operator on $w_0$-densities with
a given  ``non-singular'' weight  $w_0\neq 0,\frac{1}{2}, 1$. A
pencil can be recovered from an operator on functions up to
$w(w-1)f$, where $f$ is a scalar.

Consider odd operators. To them correspond odd brackets. Notice
that an odd symmetric bracket is transformed into an even
antisymmetric bracket by the parity shift. If $\D$ is odd  and
$\ord \D\leq 2$, then $\ord \D^2\leq 3$. The condition $\ord
\D^2\leq 2$ is equivalent to the Jacobi identity for the bracket
generated by $\D$. In this case $D=(S,\_)$ is a differential in
$\fun(T^*M)$, and for an upper connection  $\g=\g^a p_a$ the
notion of curvature makes sense: $F=D\g$.

\begin{thm}
An operator $\D$ is a derivation of the generated bracket
(equivalently, $\ord \D^2\leq 1$) if and only if $D\g=0$.
\end{thm}

The flatness of $\g$  is ``subsumed'' by the Jacobi identity for a
bracket in  $\V(M)$:

\begin{thm}
The Jacobi identity for an odd bracket in  $\V(M)$ is equivalent
to the equations $(S,S)=0$, $(S,\g)=0$, $(S,\lt)+(\g,\g)=0$,
$(\g,\lt)=0$. (In such case also $\D^2_w=\lie{X}$  for some $X\in
\Vect (M)$, which is automatically Poisson.)
\end{thm}

\begin{coro}
If the odd bracket on $M$ specified by  $S^{ab}$ is
non-degenerate, then the Jacobi identity for the bracket in
$\V(M)$ implies $\g^a=S^{ab}\g_b$, $\lt=\g^a\g_a$, where
$\g_a=-\p_a\ln \rho$, and $\rho=e^{\Act}$ is a local volume form.
Hence, $\D$ is defined by a bracket on  $M$ and an ``effective
action'' $\Act$.
\end{coro}

The condition $\D^2=0$ in the odd symplectic case gives the
Batalin--Vilkovisky equation for an action $\Act$. In the general
odd Poisson case the ``Batalin--Vilkovisky equations'' are written
for a pair of flat connections and possess a groupoid property.
They describe a change of the operator  $\Delta^2$ on functions or
$\Delta$ on half-densities~\cite{tv:laplace1,tv:laplace2}.

The theory can be generalized for operators and brackets of
nonzero weight. It is interesting to consider a generalization for
operators of higher order   (where homotopy algebras should
appear).

\def\cprime{$'$}

\end{document}